\documentclass[10pt,letterpaper,fleqn,oneside]{article}
\pdfoutput=1

\usepackage{sectsty}
\usepackage{caption}
\usepackage{times}
\usepackage{amssymb,amsfonts,amsmath,amscd,amsthm}
\usepackage[pdftex,colorlinks]{hyperref}
\usepackage[pdftex]{graphicx}
\usepackage{fnpos}
\usepackage{subfigure}
\usepackage{verbatim}
\usepackage{fancyhdr}
\usepackage[numbers]{natbib}

\usepackage{url}

\usepackage[usenames,dvipsnames]{color}
\usepackage{soul}   

\newcommand{\bbm}{\begin{bmatrix}}
\newcommand{\ebm}{\end{bmatrix}}
\newcommand{\norm}[2]{\left|\left| #1 \right|\right|_{#2}}
\newcommand{\expect}[1]{E\left[#1\right]}
\newcommand{\prob}[1]{P\left[#1\right]}

\newcommand{\dia}[0]{d}
\newcommand{\dmin}[0]{\dia_{\rm min}}
\newcommand{\stateh}[0]{\boldsymbol{\mu}}

\makeFNbottom \makeFNbelow

\newcommand{\UTIASprogram}{N/A}
\newcommand{\UTIAStitle}{On Recursive Random Prolate Hyperspheroids}
\newcommand{\UTIASdocument}{TR-2014-JDG002}
\newcommand{\UTIASrevision}{Rev: 2.0}
\newcommand{\UTIASauthor}{J.\ D.\ Gammell}

\hypersetup{%
    pdftitle={\UTIASdocument: \UTIAStitle},
    pdfauthor={\UTIASauthor},
    pdfkeywords={},
    pdfsubject={\UTIASprogram},
    pdfstartview=FitH,%
    bookmarksopen=true,%
    breaklinks=true,%
    colorlinks=true,%
    linkcolor=blue,anchorcolor=blue,%
    citecolor=blue,filecolor=blue,%
    menucolor=blue,
    urlcolor=blue
}%

\title{\sf\bfseries \UTIAStitle }

\author{Jonathan D.\ Gammell\\
       Institute for Aerospace Studies \\
       University of Toronto\\
       4925 Dufferin Street \\
       Toronto, Ontario \\
       Canada M3H 5T6 \\
       \texttt{<jon.gammell@utoronto.ca>}
       \and
       Siddhartha S.\ Srinivasa \\
       The Robotics Institute \\
       Carnegie Mellon University\\
       5000 Forbes Avenue\\
       Pittsburgh, Pennsylvania\\
       USA 15213-3890 \\
       \texttt{<siddh@cs.cmu.edu>}
       \and
       Timothy D.\ Barfoot \\
       Institute for Aerospace Studies \\
       University of Toronto\\
       4925 Dufferin Street \\
       Toronto, Ontario \\
       Canada M3H 5T6 \\
       \texttt{<tim.barfoot@utoronto.ca>}}

\date{}

\graphicspath{{figs/}}

\begin{document}

\fancypagestyle{plain}{%
    \fancyhf{}%
    \fancyfoot[C]{}%
    \fancyhead[R]{\begin{tabular}[b]{r}\small\sf \UTIASdocument\\
         \small\sf\UTIASrevision\\
         \small\sf\today \end{tabular}}%
    \renewcommand{\headrulewidth}{0pt}
    \renewcommand{\footrulewidth}{0pt}%
}

\pagestyle{fancy}
\allsectionsfont{\sf\bfseries}
\renewcommand{\captionlabelfont}{\sf\bfseries}
\reversemarginpar
\renewcommand{\labelitemi}{--}
\renewcommand{\labelitemii}{--}
\renewcommand{\labelitemiii}{--}


\rhead{ \begin{tabular}[b]{r}\small\sf \UTIASdocument\\
        \small\sf\UTIASrevision\\
        \small\sf\today \end{tabular}}
\chead{}
\lfoot{}
\cfoot{\thepage}
\rfoot{}
\renewcommand{\headrulewidth}{0pt}
\renewcommand{\footrulewidth}{0pt}


\maketitle

\begin{abstract}
This technical note analyzes the properties of a random sequence of prolate hyperspheroids with common foci. Each prolate hyperspheroid in the sequence is defined by a sample drawn randomly from the previous volume such that the sample lies on the new surface (Fig.~\ref{fig:randEllipse}). Section \ref{sec:coords} defines the prolate hyperspheroid coordinate system and the resulting differential volume, Section \ref{sec:expect} calculates the expected value of the new transverse diameter given a uniform distribution over the existing prolate hyperspheroid, and Section \ref{sec:converge} calculates the convergence rate of this sequence. For clarity, the differential volume and some of the identities used in the integration are verified in Appendix \ref{appx:volume} through a calculation of the volume of a general prolate hyperspheroid.
\end{abstract}

\begin{figure}[tb]
	\centering
	\includegraphics[width=0.5\textwidth]{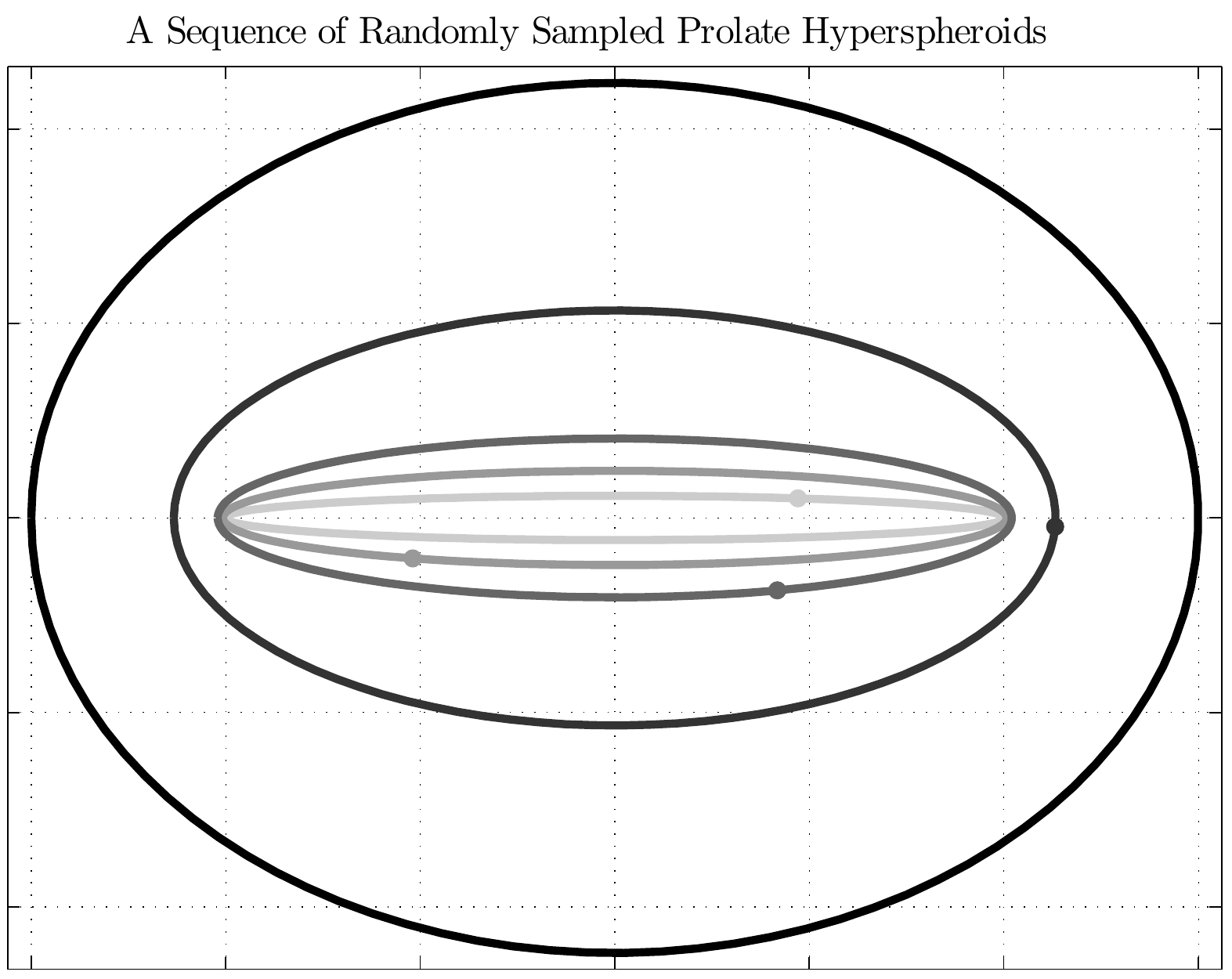}
	\caption{A series of 4 random prolate hyperspheroids generated by sampling a point from within the volume of the previous hyperspheroid. The generating sample is plotted in the same shade of grey as the resulting surface.}
	\label{fig:randEllipse}
\end{figure}
\section{Prolate Hyperspheroid Coordinate System}\label{sec:coords}
\begin{figure}[tb]
	\centering
	\includegraphics[width=0.5\textwidth]{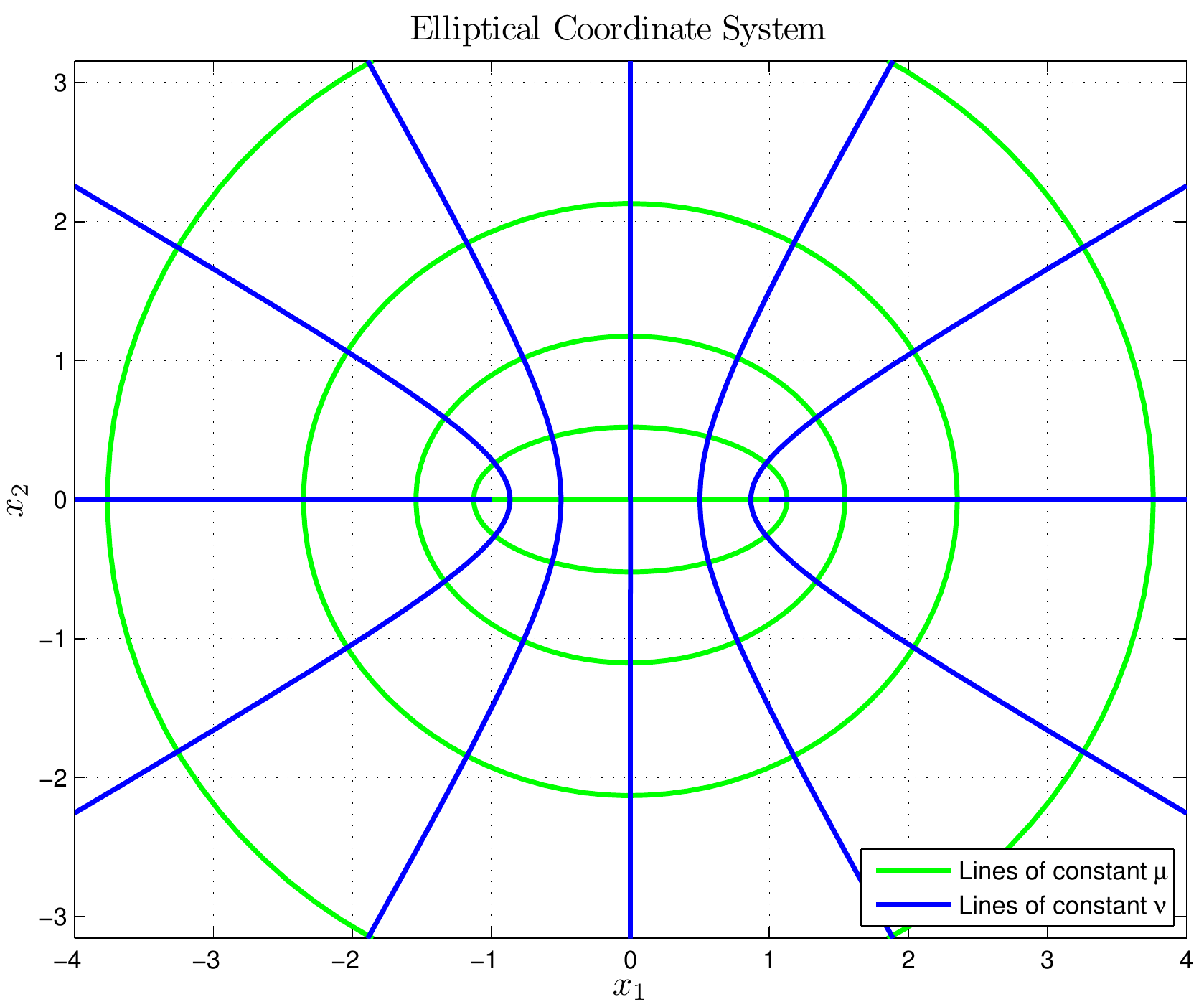}
	\caption{The 2D elliptical coordinate system for $a = 1$, note that in the 2D case, $\mu \in \left[0,\infty\right)$ and $\nu \in \left[0, 2\pi\right)$. Starting from the line on the positive $x_1$ axis and proceeding counterclockwise, the blue lines of constant $\nu$ correspond to angles of $0, \frac{\pi}{6}, \frac{2\pi}{6}, \ldots \frac{11\pi}{6}$ radians respectively. Proceeding outwards, the green lines of constant $\mu$ correspond to values of $0, 0.5, \ldots, 2$.}
	\label{fig:ellipticalCoords}
\end{figure}

Let $\left( x_1, x_2, \ldots, x_n \right)$ be the Cartesian coordinates of an $\mathbb{R}^n$ coordinate system, then we can define a prolate hyperspheroid coordinate system, $\left(\mu, \nu, \psi_1, \psi_2, \ldots, \psi_{n-2} \right)$, parameterized on $a$ as
\begin{align}\label{eqn:cart}
    x_1 &=: a\cosh\mu\cos\nu,\\
    x_2 &=: a\sinh\mu\sin\nu\cos\psi_1,\nonumber\\
    x_3 &=: a\sinh\mu\sin\nu\sin\psi_1\cos\psi_2,\nonumber\\
    &\vdots\nonumber\\
    x_{n-1} &=: a\sinh\mu\sin\nu\sin\psi_1\sin\psi_2\ldots\sin\psi_{n-3}\cos\psi_{n-2},\nonumber\\
    x_{n} &=: a\sinh\mu\sin\nu\sin\psi_1\sin\psi_2\ldots\sin\psi_{n-3}\sin\psi_{n-2},\nonumber
\end{align}
where the foci of the prolate hyperspheroids occur in Cartesian coorinates at $\left(\pm a, 0, \ldots, 0\right)$ and the transverse diameter of the prolate hyperspheroid on which a given point lies, $\dia$, is of length
\begin{align}\label{eqn:dia}
    \dia = \dmin\cosh\mu,
\end{align}
with the minimum transverse diameter, $\dmin$, defined as the distance between the foci, $\dmin = 2a$.
These coordinates can be viewed as a $2D$ elliptical coordinate system, $\mu, \nu$, rotated by spherical coordinates, $\psi_1, \psi_2,\ldots,\psi_{n-2}$. The coordinates take the values $\mu \in \left[0, \infty \right)$, $\nu \in \left[0,\pi \right]$, $\psi_1, \psi_2, \ldots, \psi_{n-3} \in \left[0, \pi \right]$, and $\psi_{n-2} \in \left[0, 2\pi\right)$, except in the 2D case when $\nu \in \left[0,2\pi\right)$ (Fig.~\ref{fig:ellipticalCoords}).

The basis vectors of this curvilinear coordinate system, $\mathbf{e}_{\mu}, \mathbf{e}_{\nu}, \mathbf{e}_{\psi_1}, \ldots, \mathbf{e}_{\psi_{n-2}} $, are defined as the partial derivatives of \eqref{eqn:cart} with respect to the respective prolate hyperspheroid coordinate,
\begin{align*}
    \mathbf{e}_{\mu} &:= \bbm \frac{\partial x_1}{\partial\mu} & \frac{\partial x_2}{\partial\mu} & \ldots & \frac{\partial x_n}{\partial\mu}\ebm^T\mbox{, etc.,}
\end{align*}
from which we can calculate the scale factors, $h_{\mu}, h_{\nu}, h_{\psi_1}, \ldots, h_{\psi_{n-2}}$, as
\begin{align}\label{eqn:scaleDefn}
    h_{\mu} := \norm{\mathbf{e}_{\mu}}{2}\mbox{, etc.}
\end{align}
The differential unit of volume, $dV$, is then
\begin{align}\label{eqn:dVdef}
    dV := h_{\mu}d\mu\, h_{\nu}d\nu\, h_{\psi_1}d\psi_1\, \ldots h_{\psi_{n-2}}d\psi_{n-2} = h_{\mu}h_{\nu} d\mu\, d\nu\, \prod_{i=1}^{n-2}\left(h_{\psi_i}d\psi_i\right).
\end{align}

\subsection{First Scale Factor}
The partial derivatives of \eqref{eqn:cart} with respect to $\mu$ are
\begin{align*}
    \frac{\partial x_1}{\partial\mu} &= -a\sinh\mu\cos\nu,\\
    \frac{\partial x_2}{\partial\mu} &= a\cosh\mu\sin\nu\cos\psi_1,\\
    \frac{\partial x_3}{\partial\mu} &= a\cosh\mu\sin\nu\sin\psi_1\cos\psi_2,\\
    &\vdots\\
    \frac{\partial x_{n-1}}{\partial\mu} &= a\cosh\mu\sin\nu\sin\psi_1\ldots\sin\psi_{n-3}\cos\psi_{n-2},\\
    \frac{\partial x_{n}}{\partial\mu} &= a\cosh\mu\sin\nu\sin\psi_1\ldots\sin\psi_{n-3}\sin\psi_{n-2}.
\end{align*}
From \eqref{eqn:scaleDefn} the scale factor, $h_\mu$, is then
\begin{align*}
    h_\mu = a\left(\vphantom{\cosh^2\mu}\right.& \sinh^2\mu\cos^2\nu + \cosh^2\mu\sin^2\nu\cos^2\psi_1 + \cosh^2\mu\sin^2\nu\sin^2\psi_1\cos^2\psi_2 + \ldots\\
             &\left.+ \cosh^2\mu\sin^2\nu\sin^2\psi_1\ldots\sin^2\psi_{n-3}\cos^2\psi_{n-2} + \cosh^2\mu\sin^2\nu\sin^2\psi_1\ldots\sin^2\psi_{n-3}\sin^2\psi_{n-2}\right)^{1/2}.
\end{align*}
Using the Pythagorean trigonometric identity and the hyperbolic analogue,
\begin{align}\label{eqn:sqrcos}
    \sin^2b + \cos^2b &\equiv 1\\
    \cosh^2a - \sinh^2a &\equiv 1\nonumber,
\end{align}
for every instance of $\cos^2\left(\cdot\right)$ and $\cosh^2\left(\cdot\right)$ yields
\begin{align*}
    h_\mu = a\left(\vphantom{\sinh^\mu}\right.
        & \overbrace{\sinh^2\mu - \sinh^2\mu\sin^2\nu}^{x_1 \mbox{ term}}
        + \overbrace{\sin^2\nu + \sinh^2\mu\sin^2\nu - \sin^2\nu\sin^2\psi_1 - \sinh^2\mu\sin^2\nu\sin^2\psi_1}^{x_2 \mbox{ term}}\\
        &+ \overbrace{\sin^2\nu\sin^2\psi_1 + \sinh^2\mu\sin^2\nu\sin^2\psi_1 - \sin^2\nu\sin^2\psi_1\sin^2\psi_2 - \sinh^2\mu\sin^2\nu\sin^2\psi_1\sin^2\psi_2}^{x_3 \mbox{ term}}
        + \ldots\\
        &+ \overbrace{\sin^2\nu\sin^2\psi_1\ldots\sin^2\psi_{n-3} + \sinh^2\mu\sin^2\nu\sin^2\psi_1\ldots\sin^2\psi_{n-3}}^{x_{n-1} \mbox{ term}}\\
        & \overbrace{-\sin^2\nu\sin^2\psi_1\ldots\sin^2\psi_{n-3}\sin^2\psi_{n-2} - \sinh^2\mu\sin^2\nu\sin^2\psi_1\ldots\sin^2\psi_{n-3}\sin^2\psi_{n-2}}^{x_{n-1} \mbox{ term}}\\
        &+ \overbrace{\sin^2\nu\sin^2\psi_1\ldots\sin^2\psi_{n-3}\sin^2\psi_{n-2} + \sinh^2\mu\sin^2\nu\sin^2\psi_1\ldots\sin^2\psi_{n-3}\sin^2\psi_{n-2}}^{x_{n} \mbox{ term}}\left. \vphantom{\sinh^\mu} \right)^{\frac{1}{2}},
\end{align*}
about which we make the following observations:
\begin{itemize}
\item The second term of the $x_1$ grouping cancels with the second term of the $x_2$ grouping.
\item For the intermediate terms, $x_2,x_3,\ldots x_{n-1}$, the third and fourth term of each $x_j$ grouping cancels out with the first and second terms of the subsequent $x_{j+1}$ grouping.
\item The two terms of the final $x_n$ grouping cancel with the third and fourth terms of the $x_{n-1}$ grouping.
\end{itemize}
These observations allow us to finally write
\begin{align}\label{eqn:hmu}
h_\mu = a\left(\sinh^2\mu + \sin^2\nu \right)^\frac{1}{2}.
\end{align}

\subsection{Second Scale Factor}
The partial derivatives of \eqref{eqn:cart} with respect to $\nu$ are
\begin{align*}
    \frac{\partial x_1}{\partial\nu} &= -a\cosh\mu\sin\nu,\\
    \frac{\partial x_2}{\partial\nu} &= a\sinh\mu\cos\nu\cos\psi_1,\\
    \frac{\partial x_3}{\partial\nu} &= a\sinh\mu\cos\nu\sin\psi_1\cos\psi_2,\\
    &\vdots\\
    \frac{\partial x_{n-1}}{\partial\nu} &= a\sinh\mu\cos\nu\sin\psi_1\ldots\sin\psi_{n-3}\cos\psi_{n-2},\\
    \frac{\partial x_{n}}{\partial\nu} &= a\sinh\mu\cos\nu\sin\psi_1\ldots\sin\psi_{n-3}\sin\psi_{n-2}.
\end{align*}
From \eqref{eqn:scaleDefn} the scale factor, $h_\nu$, is then
\begin{align*}
    h_\nu = a\left(\vphantom{\cosh^2\mu}\right.& \cosh^2\mu\sin^2\nu + \sinh^2\mu\cos^2\nu\cos^2\psi_1 + \sinh^2\mu\cos^2\nu\sin^2\psi_1\cos^2\psi_2 + \ldots\\
             &\left.+ \sinh^2\mu\cos^2\nu\sin^2\psi_1\ldots\sin^2\psi_{n-3}\cos^2\psi_{n-2} + \sinh^2\mu\cos^2\nu\sin^2\psi_1\ldots\sin^2\psi_{n-3}\sin^2\psi_{n-2}\right)^{\frac{1}{2}}.
\end{align*}
Once again making use of the identities \eqref{eqn:sqrcos} for every instance of $\cos^2\left(\cdot\right)$ and $\cosh^2\left(\cdot\right)$ gives
\begin{align*}
    h_\nu = a\left(\vphantom{\sinh^2\nu}\right.
        &\overbrace{\sin^2\nu + \sinh^2\mu\sin^2\nu}^{x_1 \mbox{ term}}
        + \overbrace{\sinh^2\mu - \sinh^2\mu\sin^2\nu - \sinh^2\mu\sin^2\psi_1 + \sinh^2\mu\sin^2\nu\sin^2\psi_1}^{x_2 \mbox{ term}}\\
        &+ \overbrace{\sinh^2\mu\sin^2\psi_1 - \sinh^2\mu\sin^2\nu\sin^2\psi_1 - \sinh^2\mu\sin^2\psi_1\sin^2\psi_2 + \sinh^2\mu\sin^2\nu\sin^2\psi_1\sin^2\psi_2}^{x_3 \mbox{term}}
        + \ldots\\
        &+ \overbrace{\sinh^2\mu\sin^2\psi_1\ldots\sin^2\psi_{n-3} - \sinh^2\mu\sin^2\nu\sin^2\psi_1\ldots\sin^2\psi_{n-3} }^{x_{n-1}\mbox{ term}}\\
        & \overbrace{-\sinh^2\mu\sin^2\psi_1\ldots\sin^2\psi_{n-3}\sin^2\psi_{n-2} + \sinh^2\mu\sin^2\nu\sin^2\psi_1\ldots\sin^2\psi_{n-3}\sin^2\psi_{n-2}}^{x_{n-1}\mbox{ term}}\\
        &+ \overbrace{\sinh^2\mu\sin^2\psi_1\ldots\sin^2\psi_{n-3}\sin^2\psi_{n-2} - \sinh^2\mu\sin^2\nu\sin^2\psi_1\ldots\sin^2\psi_{n-3}\sin^2\psi_{n-2}}^{x_{n}\mbox{ term}}\left.\vphantom{\sinh^2\nu}\right)^\frac{1}{2},
\end{align*}
about which we make the following observations:
\begin{itemize}
\item The second term of the $x_1$ grouping cancels with the second term of the $x_2$ grouping.
\item For the intermediate terms, $x_2,x_3,\ldots x_{n-1}$, the third and fourth term of each $x_j$ grouping cancels out with the first and second terms of the subsequent $x_{j+1}$ grouping.
\item The two terms of the final $x_n$ grouping cancel with the third and fourth terms of the $x_{n-1}$ grouping.
\end{itemize}
These observations allow us to finally write
\begin{align}\label{eqn:hnu}
h_\nu = a\left(\sinh^2\mu + \sin^2\nu \right)^\frac{1}{2}.
\end{align}

\subsection{Intermediate Scale Factors}
Noting that the dependence of the Cartesian coordinates \eqref{eqn:cart} on the intermediate terms, $\psi_1, \psi_2, \ldots, \psi_{n-3}$, follow a common form, we can write a general expression for their derivatives with respect to $\psi_i$ as
\begin{align*}
    \frac{\partial x_1}{\partial\psi_i} &= \frac{\partial x_2}{\partial\psi_i} = \ldots \frac{\partial x_i}{\partial\psi_i} = 0,\\
    \frac{\partial x_{i+1}}{\partial \psi_i} &= -a\sinh\mu\sin\nu\sin\psi_1\ldots\sin\psi_{i-1}\sin\psi_i, \\
    \frac{\partial x_{i+2}}{\partial \psi_i} &= a\sinh\mu\sin\nu\sin\psi_1\ldots\sin\psi_{i-1}\cos\psi_i\cos\psi_{i+1}, \\
    \frac{\partial x_{i+3}}{\partial \psi_i} &= a\sinh\mu\sin\nu\sin\psi_1\ldots\sin\psi_{i-1}\cos\psi_i\sin\psi_{i+1}\cos\psi_{i+2}, \\
    \vdots\\
    \frac{\partial x_{n-1}}{\partial \psi_i} &= a\sinh\mu\sin\nu\sin\psi_1\ldots\sin\psi_{i-1}\cos\psi_i\sin\psi_{i+1}\ldots\sin\psi_{n-3}\cos\psi_{n-2}, \\
    \frac{\partial x_{n}}{\partial \psi_i} &= a\sinh\mu\sin\nu\sin\psi_1\ldots\sin\psi_{i-1}\cos\psi_i\sin\psi_{i+1}\ldots\sin\psi_{n-3}\sin\psi_{n-2}.
\end{align*}
From \eqref{eqn:scaleDefn} the scale factor, $h_{\psi_i}$, is then
\begin{align*}
    h_{\psi_i} = a\sinh\mu\sin\nu&\sin\psi_1\ldots\sin\psi_{i-1}\left(\vphantom{\sin^2\psi_{i+n}}\right. \sin^2\psi_i + \cos^2\psi_i\cos^2\psi_{i+1} + \cos^2\psi_i\sin^2\psi_{i+1}\cos^2\psi_{i+2}
                + \ldots\\
         &\left.+ \cos^2\psi_i\sin^2\psi_{i+1}\ldots\sin^2\psi_{n-3}\cos^2\psi_{n-2} + \cos^2\psi_i\sin^2\psi_{i+1}\ldots\sin^2\psi_{n-3}\sin^2\psi_{n-2}\right)^{\frac{1}{2}}.
\end{align*}
Returning to the Pythagorean trigonometric identity \eqref{eqn:sqrcos} for every instance of $\cos^2\left(\cdot\right)$ gives
\begin{align*}
    h_{\psi_i} = a\sinh\mu\sin\nu&\sin\psi_1\ldots\sin\psi_{i-1}\left(\vphantom{\sin^2\psi_{i+n}}\right.
        \overbrace{\sin^2\psi_i}^{x_{i+1}\mbox{ term}} + 
        \overbrace{1 - \sin^2\psi_i-\sin^2\psi_{i+1} + \sin^2\psi_i\sin^2\psi_{i+1}}^{x_{i+2}\mbox{ term}}\\
        &+ \overbrace{\sin^2\psi_{i+1} - \sin^2\psi_{i}\sin^2\psi_{i+1} - \sin^2\psi_{i+1}\sin^2\psi_{i+2} + \sin^2\psi_{i}\sin^2\psi_{i+1}\sin^2\psi_{i+2} }^{x_{i+3}\mbox{ term}}
        +\ldots\\
        &+ \overbrace{\sin^2\psi_{i+1}\ldots\sin^2\psi_{n-3} - \sin^2\psi_{i}\sin^2\psi_{i+1}\ldots\sin^2\psi_{n-3}}^{x_{n-1}\mbox{ term}}\\
        & \overbrace{-\sin^2\psi_{i+1}\ldots\sin^2\psi_{n-3}\sin^2\psi_{n-2} + \sin^2\psi_{i}\sin^2\psi_{i+1}\ldots\sin^2\psi_{n-3}\sin^2\psi_{n-2}}^{x_{n-1}\mbox{ term}}\\
        &+ \overbrace{\sin^2\psi_{i+1}\ldots\sin^2\psi_{n-3}\sin^2\psi_{n-2} - \sin^2\psi_{i}\sin^2\psi_{i+1}\ldots\sin^2\psi_{n-3}\sin^2\psi_{n-2}}^{x_{n}\mbox{ term}}\left.\vphantom{\sin^2\psi_{i+n}}\right)^\frac{1}{2},
\end{align*}
about which we make the following observations:
\begin{itemize}
\item The only term of the $x_i$ grouping cancels with the second term of the $x_{i+1}$ grouping.
\item For the intermediate terms, $x_{i+2},x_{i+3},\ldots x_{n-1}$, the third and fourth term of each $x_j$ grouping cancels out with the first and second terms of the subsequent $x_{j+1}$ grouping.
\item The two terms of the final $x_n$ grouping cancel with the third and fourth terms of the $x_{n-1}$ grouping.
\end{itemize}
This leaves unity inside the square root, giving
\begin{align}\label{eqn:hpsii_1}
h_{\psi_i} = a\sinh\mu\sin\nu&\sin\psi_1\ldots\sin\psi_{i-1} \quad 1 \leq i \leq n-3.
\end{align}

\subsection{Final Scale Factor}
The partial derivatives of \eqref{eqn:cart} with respect to the last prolate hyperspheroid coordinate, $\psi_{n-2}$, are
\begin{align*}
    \frac{\partial x_1}{\partial\psi_{n-2}} &= \frac{\partial x_2}{\partial\psi_{n-2}} = \ldots \frac{\partial x_{n-2}}{\partial\psi_{n-2}} = 0,\\
    \frac{\partial x_{n-1}}{\partial \psi_{n-2}} &= -a\sinh\mu\sin\nu\sin\psi_1\ldots\ldots\sin\psi_{n-3}\sin\psi_{n-2}, \\
    \frac{\partial x_{n}}{\partial \psi_{n-2}} &= a\sinh\mu\sin\nu\sin\psi_1\ldots\ldots\sin\psi_{n-3}\cos\psi_{n-2}.
\end{align*}
From \eqref{eqn:scaleDefn} the scale factor, $h_{\psi_{n-2}}$, is then
\begin{align*}
    h_{\psi_{n-2}} = a\sinh\mu\sin\nu&\sin\psi_1\ldots\sin\psi_{n-3}\left(\sin^2\psi_{n-2} + \cos^2\psi_{n-2}\right)^{\frac{1}{2}},
\end{align*}
to which we return one last time to the Pythagorean trigonometric identity \eqref{eqn:sqrcos} to get
\begin{align*}
    h_{\psi_{n-2}} = a\sinh\mu\sin\nu&\sin\psi_1\ldots\sin\psi_{n-3}.
\end{align*}
Combined with \eqref{eqn:hpsii_1}, we can now write a single expression for all $\psi_i$:
\begin{align}\label{eqn:hpsii_2}
h_{\psi_i} = a\sinh\mu\sin\nu&\sin\psi_1\ldots\sin\psi_{i-1} \quad 1 \leq i \leq n-2.
\end{align}

\subsection{Differential Volume}
Substituting \eqref{eqn:hmu}, \eqref{eqn:hnu}, and \eqref{eqn:hpsii_2} into \eqref{eqn:dVdef}, gives a final expression for the differential volume
\begin{align}\label{eqn:dV}
    dV = a^n\left(\sinh^2\mu + \sin^2\nu \right)\sinh^{n-2}\mu\sin^{n-2}\nu\sin^{n-3}\psi_1\sin^{n-4}\psi_2\ldots\sin\psi_{n-3}\,d\mu\,d\nu\,d\psi_1\,d\psi_2\,\ldots\,d\psi_{n-3}\,d\psi_{n-2}
\end{align}

\section{Expectation of the Transverse Diameter}\label{sec:expect}
Given common foci, we can calculate the expected diameter of a new prolate hyperspheroid, $\dia_{i+1}$, constrained to pass through a sample drawn from a uniform distribution over the volume, $V_i$, of the current prolate hyperspheroid with transverse diameter, $\dia_{i}$, as
\begin{align*}
    \expect{\dia_{i+1}} = \int_{V_i}\dia\left(\mu\right)f\left(\stateh\right)\,dV.
\end{align*}
The transverse diameter on which the a sample in prolate hyperspheroid coordinates lies, $\dia\left(\cdot\right)$, is given by \eqref{eqn:dia}, $dV$ is the differential volume in prolate hyperspheroid coordinates \eqref{eqn:dV}, and $f\left(\cdot\right) = \frac{1}{V_i}$ is the probability density function, with $V_i$ being the volume of the containing prolate hyperspheroid,
\begin{align}\label{eqn:Vphs}
    V_i = \zeta_n\frac{\dia_i\left(\dia_i^2-\dmin^2\right)^\frac{n-1}{2}}{2^n},
\end{align}
with $\zeta_n$ as the volume of a unit $n$-ball. Making these substitutions and rearranging the independent integrals gives
\begin{align*}
    \expect{\dia_{i+1}} = \frac{\dmin^{n+1}}{2^nV_i} \int_0^{\mu'} \int_0^\pi &\cosh\mu\left(\sinh^2\mu + \sin^2\nu \right)\sinh^{n-2}\mu\sin^{n-2}\nu\,d\mu\,d\nu\\
    &\underbrace{\int_0^\pi \sin^{n-3}\psi_1\,d\psi_1,
    \int_0^\pi \sin^{n-4}\psi_2\,d\psi_2
    \ldots
    \int_0^\pi \sin\psi_{n-3}d\psi_{n-3}
     \int_0^{2\pi} d\psi_{n-2}}_{\left(n-1\right)\zeta_{n-1}},
\end{align*}
where we have recognized that in spherical coordinates, the volume of a general $n$-ball is given by \cite{derise_ijmest92}
\begin{align}\label{eqn:zetaInt}
    V_{\rm n-ball} = \int_0^r r^{n-1}\,dr \int_0^\pi \sin^{n-2}\phi_1\,d\phi_1 \int_0^\pi \sin^{n-3}\phi_2\,d\phi_2 \ldots \int_0^\pi \sin\phi_{n-2}\,d\phi_{n-2} \int_0^{2\pi}d\phi_{n-1}.
\end{align}
This leaves us with the more manageable equation
\begin{align}\label{eqn:expectTemp1}
        \expect{\dia_{i+1}} = \frac{\left(n-1\right)\dmin^{n+1}\zeta_{n-1}}{2^nV_i} \left(\int_0^{\mu'} \cosh\mu\sinh^n\mu\,d\mu \int_0^\pi \sin^{n-2}\nu\,d\nu + \int_0^{\mu'} \cosh\mu\sinh^{n-2}\mu\,d\mu \int_0^\pi \sin^n\nu\,d\nu\right),
\end{align}
that we can integrate with the help of beta functions.

Integrals of the product of $\sin\left(\theta\right)$ and $\cos\left(\theta\right)$ over the interval $\left[0,\pi\right]$ can be expressed in terms of the beta function \cite{schaums}, $\mathrm{B}\left(\cdot,\cdot\right)$, as
\begin{align}\label{eqn:betaInt}
    \int_0^\pi \sin^{2m-1}\theta \cos^{2n-1}\theta \,d\theta \equiv B\left( m,n \right),
\end{align}
making \eqref{eqn:expectTemp1}
\begin{align}\label{eqn:expectTemp2}
    \expect{\dia_{i+1}} = \frac{\left(n-1\right)\dmin^{n+1}\zeta_{n-1}}{2^nV_i} \left( B\left(\tfrac{n-1}{2},\tfrac{1}{2}\right)\int_0^{\mu'}\cosh\mu\sinh^n\mu\,d\mu + B\left(\tfrac{n+1}{2},\tfrac{1}{2}\right)\int_0^{\mu'} \cosh\mu\sinh^{n-2}\mu\,d\mu\right).
\end{align}
Making use of the relation between the beta function and gamma function \cite{schaums}, $\Gamma\left(\cdot\right)$,
\begin{align*}
    B\left(m,n\right) \equiv \frac{\Gamma\left(m\right)\Gamma\left(n\right)}{\Gamma\left(m+n\right)}
\end{align*}
and a common identity of the gamma function \cite{schaums},
\begin{align*}
    \Gamma\left(n+1\right) \equiv n\Gamma\left(n\right),
\end{align*}
we can write $B\left(m+1,n\right)$ in terms of $B\left(m,n\right)$ as
\begin{align}\label{eqn:recBeta}
    B\left(m+1,n\right) = \frac{\Gamma\left(m+1\right)\Gamma\left(n\right)}{\Gamma\left(m+n+1\right)} = \frac{m\Gamma\left(m\right)\Gamma\left(n\right)}{\left(m+n\right)\Gamma\left(m+n\right)} = \frac{m}{m+n}B\left(m,n\right),
\end{align}
further simplifying \eqref{eqn:expectTemp2} to
\begin{align}\label{eqn:expectTemp3}
    \expect{\dia_{i+1}} = \frac{\left(n-1\right)\dmin^{n+1}\zeta_{n-1}B\left(\tfrac{n-1}{2},\tfrac{1}{2}\right)}{2^nV_i} \left(\int_0^{\mu'} \cosh\mu\sinh^n\mu\,d\mu + \frac{n-1}{n}\int_0^{\mu'} \cosh\mu\sinh^{n-2}\mu\,d\mu\right).
\end{align}
The volume of a unit $n$-ball can be expressed in terms of the gamma function \cite{huber_amm82},
\begin{align*}
    \zeta_n \equiv \frac{\Gamma\left(\frac{1}{2}\right)^n}{\Gamma\left(\frac{n}{2}+1\right)},
\end{align*}
or as a recursive function of the volume of a $\left(n-1\right)$-ball,
\begin{align*}
    \zeta_n \equiv \frac{\Gamma\left(\frac{n+1}{2}\right)\Gamma\left(\frac{1}{2}\right)}{\Gamma\left(\frac{n}{2}+1\right)}\zeta_{n-1} = B\left( \tfrac{n+1}{2},\tfrac{1}{2} \right)\zeta_{n-1},
\end{align*}
which we can use to further rearrange using \eqref{eqn:recBeta} to
\begin{align}\label{eqn:recZeta}
    \zeta_n \equiv \frac{n-1}{n}B\left(\tfrac{n-1}{2},\tfrac{1}{2}\right)\zeta_{n-1}.
\end{align}
Substituting \eqref{eqn:recZeta} into \eqref{eqn:expectTemp3} gives
\begin{align*}
    \expect{\dia_{i+1}} = \frac{n\dmin^{n+1}\zeta_n}{2^nV_i} \left(\int_0^{\mu'} \cosh\mu\sinh^n\mu\,d\mu + \frac{n-1}{n}\int_0^{\mu'} \cosh\mu\sinh^{n-2}\mu\,d\mu\right),
\end{align*}
into which we can finally substitute \eqref{eqn:Vphs} for the volume of a prolate hyperspheroid to give
\begin{align}\label{eqn:expectInt}
    \expect{\dia_{i+1}} = \frac{n\dmin^{n+1}}{\dia_i\left(\dia_i^2 - \dmin^2\right)^\frac{n-1}{2}} \left(\int_0^{\mu'} \cosh\mu\sinh^n\mu\,d\mu + \frac{n-1}{n}\int_0^{\mu'} \cosh\mu\sinh^{n-2}\mu\,d\mu\right).
\end{align}
Given the indefinite integral \cite{schaums},
\begin{align*}
    \int \cosh x\sinh^nx\,dx = \frac{\sinh^{n+1}x}{n+1},
\end{align*}
\eqref{eqn:expectInt} becomes
\begin{align}\label{eqn:expectSinh}
    \expect{\dia_{i+1}} = \frac{n\dmin^{n+1}}{\dia_i\left(\dia_i^2 - \dmin^2\right)^\frac{n-1}{2}}  \left(\frac{\sinh^{n+1}\mu'}{n+1} + \frac{n-1}{n}\frac{\sinh^{n-1}\mu'}{n-1}\right).
\end{align}
Where $\mu'$ is given by \eqref{eqn:dia} in terms of the current transverse diameter, $\dia_i$ as,
\begin{align}\label{eqn:muPrimeCosh}
 	\cosh\mu' = \frac{\dia_i}{\dmin},
\end{align}
or if we use the identity $\cosh a = b \iff \sinh a = \sqrt{b^2 -1}$ \cite{schaums},
\begin{align}\label{eqn:muPrimeSinh}
    \sinh\mu' = \frac{1}{\dmin}\sqrt{\dia_i^2 - \dmin^2}.
\end{align}
Using \eqref{eqn:muPrimeSinh} to rearrange the remaining terms of \eqref{eqn:expectSinh} gives the final expression for the expected transverse diameter, $\expect{\dia_{i+1}}$, of a prolate hyperspheroid defined to pass through a sample taken from a uniform distribution over an earlier prolate hyperspheroid of transverse diameter $\dia_i$ with common foci as
\begin{align}\label{eqn:costExpect}
    \expect{\dia_{i+1}} = \frac{n\dia_i^2 + \dmin^2}{\left(n+1\right)\dia_i}.
\end{align}

\section{Convergence of the Expectation of the Transverse Diameter}\label{sec:converge}
The new transverse diameter, $\dia_{i+1}$ is bounded from above by the current transverse diameter, $\dia_i$,
\begin{align*}
    \dia_{i+1} \leq \dia_{i},
\end{align*}
with the diameter remaining unchanged only when the sample lies on the surface of the existing prolate hyperspheroid. As the set of such states has measure $0$, the probability of sampling a point on the surface from a uniform distribution over the volume is $0$ and as such the probability that the diameter does not decrease is also $0$,
\begin{align*}
    \prob{\dia_{i+1} = \dia_{i}} = 0.
\end{align*}
This allows us to state that the transverse diameter of the prolate hyperspheroids almost surely converges to $\dmin$,
\begin{align*}
    \prob{\lim_{i\to\infty} \dia_i = \dmin} = 1.
\end{align*}
We can then calculate the rate of the convergence, $\eta$, from \eqref{eqn:costExpect} as,
\begin{align*}
    \eta = \left. \frac{\partial \expect{\dia_{i+1}}}{\partial \dia_i} \right|_{\dia_i = \dmin},
\end{align*}
where we identified $\dmin$ as a stationary point by inspection, i.e., $\expect{\dmin} = \dmin$.
Evaluating the derivative gives
\begin{align}\label{eqn:costExpectRate}
    \eta = \frac{n-1}{n+1},
\end{align}
which as $\forall n \geq 2,\quad 0 < \frac{n-1}{n+1} < 1$, $\eta$ is always linear in convergence.

\newpage
\appendix

\section{Volume Integration Check}\label{appx:volume}
As an exercise, we should be able to recover \eqref{eqn:Vphs} from \eqref{eqn:dV} through integration,
\begin{align*}
    V = \int_VdV = \int_0^{\mu'} \int_0^\pi \int_0^\pi \int_0^\pi \ldots \int_0^\pi \int_0^{2\pi}a^n&\left(\sinh^2\mu + \sin^2\nu \right)\sinh^{n-2}\mu\sin^{n-2}\nu\\
    &\sin^{n-3}\psi_1\sin^{n-4}\psi_2\ldots\sin\psi_{n-3}\,d\mu\,d\nu\,d\psi_1\,d\psi_2\,\ldots\,d\psi_{n-3}\,d\psi_{n-2}.
\end{align*}
Rearranging the independent integrals gives
\begin{align*}
    V = \int_VdV = a^n\int_0^{\mu'} \int_0^\pi &\left(\sinh^2\mu + \sin^2\nu \right)\sinh^{n-2}\mu\sin^{n-2}\nu\,d\mu\,d\nu\\
    &\underbrace{\int_0^\pi \sin^{n-3}\psi_1\,d\psi_1,
    \int_0^\pi \sin^{n-4}\psi_2\,d\psi_2
    \ldots
    \int_0^\pi \sin\psi_{n-3}d\psi_{n-3}
     \int_0^{2\pi} d\psi_{n-2}}_{\left(n-1\right)\zeta_{n-1}},
\end{align*}
where can be simplified by the definition of the unit $n$-ball volume \eqref{eqn:zetaInt} and the beta function (\ref{eqn:betaInt}, \ref{eqn:recBeta}) to
\begin{align*}
    V = \left(n-1\right)a^n\zeta_{n-1}B\left(\tfrac{n-1}{2},\tfrac{1}{2}\right)\left( \int_0^{\mu'} \sinh^n\mu \,d\mu + \frac{n-1}{n}\int_0^{\mu'} \sinh^{n-2}\mu \,d\mu \right).
\end{align*}
Next, using the indefinite integral
\begin{align*}
    \int \sinh^n x\,dx = \frac{sinh^{n-1}x \cosh x}{n} - \frac{n-1}{n}\int\sinh^{n-2}x\,dx,
\end{align*}
we get
\begin{align*}
    V = \left(n-1\right)a^n\zeta_{n-1}B\left(\tfrac{n-1}{2},\tfrac{1}{2}\right)\left( \frac{sinh^{n-1}\mu' \cosh \mu'}{n} - \frac{n-1}{n}\int_0^{\mu'}\sinh^{n-2}\mu\,d\mu + \frac{n-1}{n}\int_0^{\mu'} \sinh^{n-2}\mu \,d\mu \right),
\end{align*}
which is simply
\begin{align*}
    V = \frac{n-1}{n}a^n\zeta_{n-1}B\left(\tfrac{n-1}{2},\tfrac{1}{2}\right) sinh^{n-1}\mu' \cosh \mu'.
\end{align*}
Recognizing the appropriate terms from the recursive definition of the unit $n$-ball \eqref{eqn:recZeta} allows us to write
\begin{align*}
    V = a^n\zeta_n \sinh^{n-1}\mu' \cosh \mu'
\end{align*}
which can be evaluated using \eqref{eqn:muPrimeCosh}, \eqref{eqn:muPrimeSinh}, and the fact that $\dmin = 2a$ to finally give
\begin{align}
    V = \zeta_n\frac{\dia\left(\dia^2 - \dmin^2\right)^{\frac{n-1}{2}}}{2^n},
\end{align}
which is exactly \eqref{eqn:Vphs}.

\newpage
\bibliographystyle{asrl}
\bibliography{TR-2014-JDG002}

\begin{thebibliography}{3}
\newcommand{\enquote}[1]{``#1''}
\providecommand{\natexlab}[1]{#1}

\bibitem[{DeRise(1992)}]{derise_ijmest92}
DeRise, G., \enquote{Some $n$-dimensional geometry,} \emph{International
  Journal of Mathematical Education in Science and Technology}, 23(3):371--379,
  1992.

\bibitem[{Huber(1982)}]{huber_amm82}
Huber, G., \enquote{Gamma Function Derivation of $n$-Sphere Volumes,} \emph{The
  American Mathematical Monthly}, 89(5):301--302, 1982.

\bibitem[{Spiegel et~al.(2012)Spiegel, Lipschutz, and Liu}]{schaums}
Spiegel, M.~R., Lipschutz, S., and Liu, J., \emph{Mathematical Handbook of
  Formulas and Tables}, Schaum's Outlines, McGraw-Hill, 4$^{th}$ edition, 2012.

\end{thebibliography}

\end{document}